\documentclass[12pt]{amsart}
\usepackage{amsmath,amsthm,amsfonts,amssymb,eucal}

\renewcommand {\a}{ \alpha }

\newcommand{\y}{\eta}

\newcommand{\vare}{\varepsilon}
\newcommand{\g}{\gamma}
\newcommand{\G}{\Gamma}

\newcommand{\vark}{\varkappa}

\renewcommand{\d}{\delta}
\newcommand{\D}{\Delta}
\newcommand{\s}{\sigma}
\renewcommand{\l}{\lambda}
\renewcommand{\L}{\Lambda}
\newcommand{\z}{\zeta}

\newcommand{\p}{\partial}
\newcommand{\om}{\omega}
\newcommand{\Om}{\Omega}

\newcommand{\R}{ \mathbb R}
\newcommand{\C}{ \mathbb C}
\newcommand{\N}{ \mathbb N}

\newcommand{\CL}{\mathcal L}

\newcommand{\CD}{\mathcal D}
\newcommand{\CE}{\mathcal E}
\newcommand{\CG}{\mathcal G}
\newcommand{\CH}{\mathcal H}
\newcommand{\CJ}{\mathcal J}

\newcommand{\CS}{\mathcal S}
\newcommand{\CW}{\mathcal W}
\newcommand{\CX}{\mathcal X}
\newcommand{\CM}{\mathcal M}

\newcommand {\GH}{\mathfrak H}

\newcommand {\GS}{\mathfrak S}
\newcommand {\GY}{\mathfrak Y}
\newcommand {\gm}{\mathfrak m}

\newcommand {\BA}{\mathbf A}

\newcommand {\BG}{\mathbf G}
\newcommand {\BH}{\mathbf H}
\newcommand {\BI}{\mathbf I}
\newcommand {\BK}{\mathbf K}

\newcommand {\BT}{\mathbf T}
\newcommand {\BV}{\mathbf V}
\newcommand {\BZ}{\mathbf Z}

\newcommand{\CA}{\mathcal A}

\newcommand{\tC}{\mathrm C}
\newcommand{\tD}{\mathrm D}
\newcommand{\tH}{\mathrm H}
\newcommand{\tY}{\mathrm Y}

\newcommand{\tL}{\mathrm L}

\newcommand{\wt}{\widetilde}
\newcommand{\wh}{\widehat}

 \DeclareMathOperator{\im}{Im}
\DeclareMathOperator{\re}{Re} 
\DeclareMathOperator{\res}{\restriction}

\DeclareMathOperator{\rank}{rank}
\DeclareMathOperator{\meas}{meas}

\newtheorem{thm}{Theorem}[section]
\newtheorem{cor}[thm]{Corollary}
\newtheorem{lem}[thm]{Lemma}
\newtheorem{prop}[thm]{Proposition}

\theoremstyle{definition}
\theoremstyle{remark}

\numberwithin{equation}{section}

\newcommand{\thmref}[1]{Theorem~\ref{#1}}

\newcommand{\bsymb}{\boldsymbol}
\newcommand{\jp}{\bsymb{\bigl[}}
\newcommand{\pj}{\bsymb{\bigl]}}
\newcommand{\nul}{Ker}

\begin{document}

\title[spectrum in a model of irreversible quantum graph]
{On the absolutely continuous spectrum in a model
 of irreversible quantum graph}
\author[S.N. Naboko]{Sergey N. Naboko}
\address{Department of Mathematical Physics\\St.Petersburg State University\\
St. Petergoff 198904 St. Peterburg\\ Russia} \email
{naboko@snoopy.phys.spbu.ru}
\author[M. Solomyak]{Michael Solomyak}
\address{Department of Mathematics\\The Weizmann Institute of Science\\
Rehovot 76100\\Israel} \email{michail.solomyak@weizmann.ac.il}
\date{25.05.2004}

\begin{abstract}
A family $\BA_\a$ of differential operators depending on a real
parameter $\a\ge 0$ is considered. This family was suggested by
Smilansky as a model of an irreversible quantum system. We find
the absolutely continuous spectrum $\s_{a.c.}$ of the operator
$\BA_\a$ and its multiplicity for all values of the parameter. The
spectrum of $\BA_0$ is purely a.c. and admits an explicit
description. It turns out that for $\a<\sqrt 2$ one has
$\s_{a.c.}(\BA_\a)= \s_{a.c.}(\BA_0)$, including the multiplicity.
For $\a\ge\sqrt2$ an additional branch of absolutely continuous
spectrum arises, its source is an auxiliary Jacobi matrix which is
related to the operator $\BA_\a$. This birth of an extra-branch of
a.c. spectrum is the exact mathematical expression of the effect
which was interpreted by Smilansky as irreversibility.
\end{abstract}
\maketitle

\section{Introduction}\label{intro}
In this paper we study the spectrum of a family $\BA_\a$ of
differential operators in the space $\tL^2(\R^2)$. Each operator
$\BA_\a$ is defined by the same differential expression
\begin{equation}\label{prel.eq}
\CA U=-U''_{xx}+\frac1{2}\bigl(-U''_{qq}+q^2 U\bigr).
\end{equation}
The parameter $\a\in\R$ appears in the ``transmission condition''
across the line $x=0$:
\begin{equation}\label{prel.tran}
U'_x(0+,q)-U'_x(0-,q)=\a q\,U(0,q),\qquad q\in\R.
\end{equation}
As we shall show in \thmref{sa.t1}, for any $\a$ the operator
$\BA_\a$ has a unique natural self-adjoint realization. The
replacement $\a\mapsto -\a$ corresponds to the change of variables
$q\mapsto-q$ which does not affect the spectrum. For this reason,
below we shall discuss only non-negative $\a$.

The family $\BA_\a$ was suggested by the physicist Smilansky in
\cite{sm} as a model of an irreversible quantum system. He carried
out a formal computation of the scattering matrix for the pair
$(\BA_0,\BA_\a)$ and showed that this matrix is unitary only if
$\a<\sqrt 2$. The loss of unitarity of the scattering matrix for
large values of $\a$ was interpreted in \cite{sm} as
irreversibility of the system.

First rigorous mathematical results on the family $\BA_\a$ were
obtained in the paper \cite{sol}, inspired by \cite{sm}. The
present paper can be considered as the second part of \cite{sol},
but it can be read independently.

The family $\BA_\a$ exhibits many unusual features, partly
revealed in \cite{sm} and \cite{sol}. The most important of them
is a phase transition at the point $\a=\sqrt 2$: the spectral
properties of $\BA_\a$ for $\a<\sqrt 2$ and for $\a>\sqrt 2$ are
quite different. In what follows we refer to the values $\a<\sqrt
2$ as ``small'' and to $\a>\sqrt 2$ as ``large''.

The spectrum of the operator $\BA_0$ can be easily described via
separation of variables. It is absolutely continuous, fills the
half-line $[1/2,\infty)$, and its multiplicity function is given
by eq. \eqref{A0.mulsac} in Section \ref{prel}. It seems natural
to study the spectrum $\s(\BA_\a)$ for $\a>0$ with the help of
the perturbation theory of quadratic forms. The standard
assumption in this type of problems is relative compactness (see
\cite{rs}) of the perturbation with respect to the unperturbed
quadratic form, i.e. to the one of the operator $\BA_0$. However,
in our case this property is violated: the perturbation is only
form-bounded but not form-compact. This was shown in \cite{sol}.

As a rule, for the form-bounded perturbations the quadratic form
approach does not give much information. Nevertheless, a rather
complete description of the essential spectrum and of the point
spectrum of $\BA_\a$ for small $\a$ was obtained in \cite{sol} and
\cite{sol1} by means of this approach. In particular, it was shown
in \cite{sol1} that for any $\a<\sqrt 2$ the spectrum  of the
operator $\BA_\a$ below the threshold $1/2$ consists of a finite
number of eigenvalues. This number grows indefinitely as
$\a\nearrow\sqrt 2$ and has regular asymptotics of the order
$O\bigl((\sqrt 2-\a)^{-1/2}\bigr)$. If $\a\ge\sqrt 2$, then the
point spectrum of $\BA_\a$ is empty (see \thmref{main.p}).
\vskip0.2cm Our goal in this paper is to study the absolutely
continuous spectrum $\s_{a.c.}$ of the operators $\BA_\a$ for all
values of the parameter $\a$. This problem was not dealt with in
\cite{sol}. A certain Jacobi operator in the space $\ell^2$ is
involved in the description of $\s_{a.c.}(\BA_\a)$, namely
\begin{gather*}
    \CJ_0(\mu):\{C_n\}\mapsto\{d_{n+1}C_{n+1}+(2n+1)\mu C_n+d_n
    C_{n-1}\},\\ d_n=n^{1/2}(n^2-1/4)^{1/4}.
\end{gather*}
Here $\mu>0$ is an auxiliary parameter. It is often convenient to
use it along with $\a$. Our main result, \thmref{main.abs}, states
that
\begin{gather}
\s_{a.c.}(\BA_\a)=\s_{a.c.}(\BA_0)\cup\s_{a.c.}(\CJ_0(\sqrt2/\a)),\label{int.1}\\
\gm_{a.c.}(E;\BA_\a)=\gm_{a.c.}(E;\BA_0)+\gm_{a.c.}(E;\CJ_0(\sqrt2/\a)),
\qquad a.e.\ E\in\R. \label{int.2}
\end{gather}
Here the symbol $\gm_{a.c.}$ stands for the multiplicity function
of the absolutely continuous spectrum.

\vskip0.2cm The spectrum of $\CJ_0(\mu)$ is discrete for $\mu>1$
and is purely absolutely continuous for $\mu\le 1$. Moreover, we
show in \thmref{res.j0spec} that
\begin{gather*}
\s(\CJ_0(\mu))=(-\infty,\infty)\ {\text{for}}\ \mu<1,\qquad
\s(\CJ_0(1))=[0,\infty);\\
\gm_{a.c.}(E;\CJ_0(\mu))=1\qquad a.e.\ {\text{on}}\
\s(\CJ_0(\mu)).
\end{gather*}
Thus, the equalities \eqref{int.1} and \eqref{int.2} give the
complete description of the absolutely continuous spectrum of the
operators $\BA_\a$. Namely, for small $\a$ it coincides with
$\s_{a.c.}(\BA_0)$, including equality of the multiplicities. For
large $\a$ a new branch of the absolutely continuous spectrum of
multiplicity $1$ adds to $\s_{a.c.}(\BA_0)$, its source is the
Jacobi matrix $\CJ_0(\mu)$. This birth of an additional branch of
the a.c. spectrum is the exact mathematical expression of the
effect which was interpreted in \cite{sm} as irreversibility.

\bigskip

The family $\BA_\a$ is a striking example of a problem which, in
spite of its seeming simplicity, exhibits many unexpected effects.
This refers to both the point spectrum and the absolutely
continuous spectrum. For this reason we believe that the detailed
analysis of this family is of general interest. Note that one
important question remains unanswered. Namely, our method does not
check whether $\s(\BA_\a)$ has singular continuous component.

\bigskip

In the course of the proof of the equalities \eqref{int.1} and
\eqref{int.2} we use the tools coming from different parts of the
spectral theory. In Section \ref{res} we obtain a convenient
representation of the operator
\begin{equation}\label{int.rr}
    (\BA_\a-\L)^{-1}-(\BA_0-\L)^{-1}.
\end{equation}
This representation involves some matrix-valued function
$\CJ_0(\L;\mu)$ which arises in a natural way when looking for
formal (i.e. not necessarily lying in $\tL^2(\R^2)$) solutions of
the equation $\BA_\a V=\L V$. The function $\CJ_0(\L;\mu)^{-1}$ is
close, in an appropriate sense, to the resolvent of the Jacobi
operator $\CJ_0(\mu)$ involved in \eqref{int.1} and \eqref{int.2}.
This allows us to find a connection between the boundary behaviour
of these two matrix-valued functions as $\L$ approaches the real
line. Technically, this is the most difficult part of the paper.
Here we make use of theory of analytic operator-valued functions.

In Section \ref{ml} we present \thmref{Aa.l1b} which relates the
a.c. spectrum of a self-adjoint operator, and also its
multiplicity function, to the jump of its "bordered resolvent"
across the real line. We could not find this result in its full
generality in the literature. For the reader's convenience, we
present its proof in Appendix \ref{ap2}.

The above mentioned representation of the operator \eqref{int.rr}
leads to an equality which expresses the jump of the bordered
resolvent of the operator $\BA_\a$ through the similar
characteristics of $\CJ_0(\mu)$. As soon as this is done,
\thmref{Aa.l1b} applies and gives the equality \eqref{int.2}. This
scheme is especially transparent for the values $E<1/2$. In order
to include the values $E>1/2$, we need an additional technical
trick.

 \vskip0.2cm Below we briefly
describe the structure of the paper. Sections \ref{prel} --
\ref{sa} contain the necessary technical material. Our main result
on the absolutely continuous spectrum $\s_{a.c.}(\BA_\a)$ for
different values of $\a$, \thmref{main.abs}, is formulated in
Section \ref{main}. Its proof is given in Sections \ref{res} --
\ref{big}.We also present \thmref{main.p} on the point spectrum of
$\BA_\a$. The latter material is mostly borrowed from the papers
\cite{sol} and \cite{sol1}.

In Section \ref{conc} we discuss the possibility to extend the
results to operators on quantum graphs. As a matter of fact, the
family $\BA_\a$ was suggested by Smilansky in \cite{sm} for this,
more general situation.

Two appendices are devoted to the proofs of \thmref{res.j0spec}
and \thmref{Aa.l1b}.

\bigskip

The notation used in the paper is mostly standard. We denote
$\N_0=\{0\}\cup\N$. The symbol $(s)\text{-}\lim$ stands for the
strong limit of operators. Abbreviation "a.e." always means almost
everywhere with respect to the Lebesgue measure. The symbols
$\tH^l$ stand for Sobolev spaces. Other necessary notations are
introduced in the course of presentation.

\section{Reduction to an infinite system of ODE}\label{prel}
Equation \eqref{prel.eq} involves the harmonic oscillator in the
variable $q$. For this reason, it is convenient to represent the
functions $U\in \GH$ as
\begin{equation*}%\label{prel.6}
U(x,q)=\sum_{n\in\N_0}u_n(x)\chi_n(q),
\end{equation*}
where $\chi_n$ are Hermite functions, normalized in $\tL^2(\R)$.
We often identify a function $U(x,q)$ with the sequence
$\{u_n(x)\}$ and write $U\sim\{u_n\}$. This identification is a
unitary mapping of the space $\tL^2(\R^2)$ onto the Hilbert space
$\GH=\ell^2\bigl(\N_0,\tL^2(\R)\bigr)$ or, equivalently, onto the
tensor product $\ell^2\otimes\tL^2(\R)$ with the natural Hilbert
space structure. For $U\sim\{u_n\}$ we have
\begin{equation}\label{prel.10}
\CA U\sim\{(L+n)u_n\},\qquad(L u)(x)=-u''(x)+u(x)/2,\ x\neq 0.
\end{equation}
The condition \eqref{prel.tran} gives
\begin{equation*}
\sum_{n\in\N_0}\chi_n(q)\bigl(u'_n(0+)-u'_n(0-)\bigr)=
\a\sum_{n\in\N_0}q\chi_n(q)u_n(0).
\end{equation*}
Taking into account the recurrence equation for the Hermite
functions:
\begin{equation*}
\sqrt{n+1}\chi_{n+1}(q)-\sqrt{2}q\chi_{n}(q)+\sqrt{n}\chi_{n-1}(q)=0,
\end{equation*}
we come to the system of matching conditions
\begin{equation}\label{prel.16}
\mu\bigl(u_n'(0+)-u_n'(0-)\bigr)=\sqrt{n+1}\,u_{n+1}(0)+\sqrt
n\,u_{n-1}(0)
\end{equation}
where $\mu=\sqrt2/\a$.

The operator $\BA_0$ admits separation of variables, which leads
to the complete description of its spectrum. Let $\BH$ stand for
the operator in $\tL^2(\R)$, defined as
\begin{equation}\label{A0.A}
\BH=-d^2/dx^2+1/2,\qquad Dom\,\BH=\tH^2(\R).
\end{equation}
One should distinguish between the self-adjoint operator $\BH$ and
the differential expression (formal differential operator) $L$
defined in \eqref{prel.10}. The operator $\BA_0$ splits into the
orthogonal sum of the operators $\BH+n$, $n\in\N_0$. An element
$U\sim\{u_n\}$ belongs to the domain $\bsymb{\CD}_0:=Dom\,\BA_0$
if and only if $u_n\in \tH^2(\R)$ for each $n\in\N_0$ and
\begin{equation*}
\sum_{n\in\N_0}\|(\BH+n)u_n\|^2<\infty.
\end{equation*}
Here and in what follows, unless otherwise explicitly stated, the
symbols $(.,.),\ \|.\|$ without indication of the space stand for
the scalar product and the norm in $\tL^2(\R)$.

The spectrum $\s(\BH)$ is absolutely continuous of multiplicity
$2$ and coincides with the half-line $[1/2,\infty)$. As a
consequence, the spectrum of $\BA_0$ is also absolutely
continuous, and
\begin{gather}
\s(\BA_0)=[1/2,\infty),\notag\\
\gm_{a.c.}(E;\BA_0)=2n\qquad {\text{for}} \ E\in(n-1/2,n+1/2),\
n\in\N. \label{A0.mulsac}
\end{gather}

\section{Jacobi matrices $\CJ(\L;\mu)$ and $\CJ_0(\mu)$}\label{aux}
\subsection{Preliminaries.}\label{pr}
Our main goal in this section is to define the Jacobi matrix
$\CJ_0(\mu)$ involved in \eqref{int.1}, \eqref{int.2}, and another
matrix $\CJ(\L;\mu)$, depending on two parameters $\L\in\C$ and
$\mu>0$. This latter matrix appears when analyzing the homogeneous
equation
\begin{equation}\label{aux.1}
\CA V=\L V
\end{equation}
under the matching condition \eqref{prel.tran}. In the
representation $V\sim\{v_n\}$ equation \eqref{aux.1} reduces to
the infinite system
\begin{gather}
-v_n''(x)+(n+1/2-\L)v_n(x)=0,\ x\neq 0; \qquad
n\in\N_0\label{aux.w}
\end{gather}
under the matching conditions \eqref{prel.16}. Our immediate task
is to describe all the formal solutions of this system. This means
the following. Introduce the linear space
\begin{equation}\label{aux.CW}
\CW(\R)=\bigl\{u\in \tC(\R): u\res\R_\pm\in \tH^2(\R_\pm)\bigr\}.
\end{equation}
Note that for any $u\in\CW(\R)$ the left-hand side in
\eqref{prel.16} is well-defined. We seek the solutions
$V\sim\{v_n\}$  such that $v_n\in\CW(\R)$ for each $n$, but not
necessarily $V\in\GH$.

Denote
\begin{gather}
\z_n:=\z_n(\L)=\sqrt{n+1/2-\L}.\label{A0.zn}
\end{gather}
We take the branch of the square root which is analytic in the
domain $\Om_n:=\C\setminus[n+1/2,\infty)$ and such that
$\z_n(\L)>0$ for $\L=\overline{\L}<n+1/2$, then
\begin{equation*}
\re \z_n(\L)>0,\  \im\z_n(\L)\cdot\im\L<0,\qquad \L\in \Om_n.
\end{equation*}
The subspace of continuous $\tL^2$-solutions of each equation
\eqref{aux.w} is one-dimensional, it is generated by the function
\begin{equation}\label{A0.yn}
\y_n(x;\L)=(n+1/2)^{1/4}e^{-\z_n(\L)|x|}.
\end{equation}
We choose such normalization of the vector-valued functions
$\y_n(.;\L)$, that each of them is analytic in $\Om_n$ and
\begin{equation}\label{A0.nnorm}
c_1(\L)\le \|\y_n(.;\L)\|^2\le c_2(\L), \qquad\forall n\in\N_0.
\end{equation}
We have
\begin{equation}\label{A0.jump}
\begin{cases}\y_n(0;\L)=(n+1/2)^{1/4},\\
\y_n'(0+;\L)-\y_n'(0-;\L)=-2(n+1/2)^{1/4}\z_n(\L).\end{cases}
\end{equation}

\vskip0.2cm

From \eqref{aux.w} we obtain $v_n(x)=C_n \y_n(x;\L)$, and by
\eqref{A0.jump} the matching conditions \eqref{prel.16} reduce to
the recurrence  system
\begin{multline*}
(n+1)^{1/2}(n+3/2)^{1/4}C_{n+1}+2\mu (n+1/2)^{1/4}\z_n(\L)C_n\\ +
n^{1/2}(n-1/2)^{1/4}C_{n-1}=0,\qquad n\in\N_0.
\end{multline*}
Taking into account our further needs, we multiply both sides of
this equality by $(n+1/2)^{1/4}$. As a result, we find that
equation \eqref{aux.1} is equivalent to the system
\begin{gather}
V\sim \{C_n \y_n(x;\L)\};\label{aux.equ}\\
d_{n+1}C_{n+1}+2\mu\, y_n(\L)C_n+d_nC_{n-1}=0\label{aux.difeq}
\end{gather}
where
\begin{equation}\label{res.entr}
d_n=n^{1/2}(n^2-1/4)^{1/4},\qquad y_n(\L)=(n+1/2)^{1/2}\z_n(\L).
\end{equation}

The Jacobi matrix $\CJ(\L;\mu)$ which corresponds to equation
\eqref{aux.difeq} is one of our main objects. In our notations we
do not distinguish between a Jacobi matrix and the operator which
it defines in $\ell^2(\N_0)$. Now we write the operator
$\CJ(\L;\mu)$ in a more convenient form. Given a number sequence
$\{\om_n\}_{n\in\N_0}$, let $\CD\{\om_n\}$ stand for the diagonal
operator in $\ell^2(\N_0)$ acting as
\begin{equation*}
\CD\{\om_n\}:\{r_0,r_1,\ldots\}\mapsto\{\om_0r_0,\om_1r_1,\ldots\}.
\end{equation*}
Let, in particular,
\begin{gather*}
\tD=\CD\{d_n\},\qquad \tY(\L)= \CD\{y_n(\L)\}.
\end{gather*}
Denoting by $\CS$ the operator of the forward shift in
$\ell^2(\N_0)$,
\begin{equation*}
\CS:\{r_0,r_1,\ldots\}\mapsto\{0,r_0,r_1,\ldots\},
\end{equation*}
we can re-write the operator  $\CJ(\L;\mu)$ as
\begin{gather}
\CJ(\L;\mu)=\tD\CS+\CS^*\tD +2\mu\tY(\L).\label{res.j}
\end{gather}

We also let
\begin{gather}
\CJ_0(\mu)=\tD\CS+\CS^*\tD +2\mu\tY_0,\qquad
\tY_0=\CD\{n+1/2\}.\label{res.j0}
\end{gather}
The operator $\CJ_0(\mu)$, defined initially on the set of all
sequences with a finite number of non-zero elements, is
essentially self-adjoint in $\ell^2$, and we denote by the same
symbol $\CJ_0(\mu)$ its unique self-adjoint extension. We do not
need the explicit description of its domain $Dom\, \CJ_0(\mu)$.
The next result describes the spectral properties of $\CJ_0(\mu)$
depending on $\mu$.
\begin{thm}\label{res.j0spec}
For $\mu>1$ the operator $\CJ_0(\mu)$ is positive definite and its
spectrum is discrete.

For $\mu\le 1$ the  spectrum of $\CJ_0(\mu)$ is purely absolutely
continuous and $\gm_{a.c.}(E;\CJ_0(\mu))=1$ a.e. on
$\s(\CJ_0(\mu))$. Moreover,
\begin{equation*}
\s(\CJ_0(1))=[0,\infty);\qquad \s(\CJ_0(\mu))=\R\ {\text{for}}\
\mu<1.
\end{equation*}
\end{thm}
The proof is given in Appendix A.

\vskip0.2cm

The difference $\CJ(\L;\mu)-(\CJ_0(\mu)-\mu\L)$ is a compact
operator. This follows from the equality
\begin{equation}\label{res.Psi}
\Psi(\L;\mu):=\CJ(\L;\mu)-(\CJ_0(\mu)-\mu\L)=\mu\CD\{\psi_n(\L)\}
\end{equation}
where
\begin{gather}
\psi_n(\L)=y_n(\L)-2(n+1/2-\L/2)\label{res.psi}
\\=-\L^2\bigl(4y_n(\L)+
4(n+1/2-\L/2)\bigr)^{-1}=O(n^{-1}).\notag
\end{gather}
Hence, the operator $\CJ(\L;\mu)$ is closed on the domain
$Dom\,\CJ(\L;\mu)=Dom\,\CJ_0(\mu)$. Moreover,
\begin{equation*}
\im\CJ(\L;\mu)=2\mu\im\tY(\L).
\end{equation*}
This implies an important property:
\begin{equation*}
\im \CJ(\L;\mu)<0\ {\text{for}}\ \L\in \C_+;\qquad \im
\CJ(\L;\mu)>0\ {\text{for}}\ \L\in \C_-.
\end{equation*}
Besides, the operator-valued function $\CJ(\L;\mu)$ forms a
holomorphic family of type (A) in the variable $\L$, see \cite{K},
\S VII.2.

It follows from \eqref{res.entr} that there is a constant $c>0$
such that
\begin{equation*}
\im y_n(-i\tau)\ge c\sqrt\tau,\qquad \forall n\in\N_0, \ \tau>1.
\end{equation*}
Hence, for $\tau>1$ we have
\begin{equation*}
\im\bigl(\CJ(-i\tau;\mu)-2i\mu c\sqrt\tau\bigr)\ge 0.
\end{equation*}
By a well known estimate for the dissipative operators, see e.g.
\cite{gokr}, Theorem IV.4.1, this implies a useful inequality
\begin{equation}\label{res.norm}
\|\CJ(-i\tau;\mu)^{-1}\|\le(2\mu c\sqrt\tau)^{-1}, \qquad \tau>1.
\end{equation}

\subsection{Birkhoff -- Adams theorem.} We base our analysis of the system
\eqref{aux.difeq}, and also the proof of \thmref{res.j0spec} on
the classical result due to Birkhoff and Adams, see e.g. the book
\cite{el}, Theorem 8.36. For the reader's convenience, we
reproduce the formulation of the part of this theorem which we
need below. It concerns the general recurrence system
\begin{equation}\label{aux.dieq}
C(n+1)+p_1(n)C(n)+p_2(n)C(n-1)=0
\end{equation}
where the functions $p_1(n)$ and $p_2(n)$ have asymptotic
expansions of the form
\begin{equation}\label{aux.dec}
p_1(n)\sim\sum_{j=0}^\infty a_jn^{-j},\qquad
p_2(n)\sim\sum_{j=0}^\infty b_jn^{-j},\ b_0\neq 0.
\end{equation}
Let $\l_\pm$ stand for the roots of the equation
\begin{equation*}
\l^2+a_0\l+b_0=0.
\end{equation*}
\begin{prop}\label{aux.bia} (a) Let $\l_+\neq\l_-$, then the system
\eqref{aux.dieq} has two linearly independent solutions
$\{C^{\pm}\}$ with the asymptotics
\begin{gather*}
C^{\pm}(n)\sim\l_\pm^n n^{d_\pm}, \qquad
d_\pm=\frac{a_1\l_\pm+b_1}{a_0\l_\pm+2b_0}.
\end{gather*}

(b) Let $\l_+=\l_-=\l$ but $2b_1\neq a_0a_1$. Then the system
\eqref{aux.dieq} has two linearly independent solutions
$\{C^{\pm}(n)\}$ with the asymptotics
\begin{equation}\label{aux.sol1}
C^{\pm}(n)\sim\l^n e^{\pm\d\sqrt n} n^\vark,\qquad
\d=2\sqrt{\frac{a_0a_1-2b_1}{2b_0}},\
\vark=\frac1{4}+\frac{b_1}{2b_0}.
\end{equation}
\end{prop}

Actually, Theorem 8.36 in \cite{el} describes the complete
asymptotic expansions of the solutions $C^{\pm}(n)$ of the system
\eqref{aux.difeq}, but we need only their leading terms. (Note
that there is an evident misprint in eq. (8.6.7) in \cite{el},
whose part is reproduced above in \eqref{aux.sol1}. In
\eqref{aux.sol1} this misprint is corrected.)\vskip0.2cm

We also need an identity for solutions of recurrence equations
with Jacobi matrices, of the type
\begin{equation}\label{aux.rec}
Q_{n+1}C_{n+1}+P_nC_n+Q_nC_{n-1}=0,\qquad n\in\N_0,
\end{equation}
with $Q_n$ real and $Q_0=0$. Namely,
\begin{equation}\label{aux.z}
\sum_{n=0}^N|C_n|^2\im
P_n=-Q_{N+1}\im\bigl(C_{N+1}\overline{C_N}\bigr), \qquad\forall
N\in\N.
\end{equation}
The proof is straightforward and we skip it.

\subsection{Solutions of the system \eqref{aux.difeq}.}\label{sys}
Here we apply Proposition \ref{aux.bia} to the system
\eqref{aux.difeq} that is actually equivalent to equation
\eqref{aux.1}. The system \eqref{aux.difeq} can be re-written in
the form \eqref{aux.dieq}, with the functions $p_1(n), p_2(n)$
admitting the asymptotic expansions of the type \eqref{aux.dec},
where in particular
 \begin{equation}\label{sys.1}
a_0=2\mu,\ a_1=-\mu(1+\L);\qquad b_0=1,\ b_1=-1.
\end{equation}
The following Lemma is a direct consequence of Proposition
\ref{aux.bia}.
\begin{lem}\label{sys.lem}
Let $\mu>0$ and $\L\in\Om_0$. Then the system \eqref{aux.difeq}
has two linearly independent solutions whose asymptotic behaviour
is given by
\begin{equation}\label{res.r3}
C^{\pm}_n\sim
\begin{cases}
(\mu+i\sqrt{1-\mu^2})^{\pm n}n^{-\frac1{2}\mp i\frac{\L\mu}
{2\sqrt{1-\mu^2}}},\qquad &\mu<1;\\
(-1)^n e^{\pm 2\sqrt{-\L n}}n^{-1/4},\qquad &\mu=1;\\
(-\mu+\sqrt{\mu^2-1})^{\pm n}n^{-\frac1{2}\pm \frac{\L\mu}
{2\sqrt{\mu^2-1}}},\qquad &\mu>1.
\end{cases}
\end{equation}
\end{lem}
We repeatedly use this Lemma in our further exposition.

\section{The self-adjoint  operator $\BA_\a,\ \a> 0$}\label{sa}
Here we describe the domain on which the operator given by
equations \eqref{prel.eq} -- \eqref{prel.tran} or equivalently, by
\eqref{prel.10} -- \eqref{prel.16}, is self-adjoint in $\GH$.
Consider a linear subset $\bsymb{\CD}_\a\subset\GH$: an element
$U\sim\{u_n\}$ belongs to $\bsymb{\CD}_\a$ if and only if each
component $u_n$ lies in $\CW(\R)$ (see \eqref{aux.CW}), the
conditions \eqref{prel.16} are satisfied, and
\begin{equation}\label{sa.1}
 \sum_{n\in\N_0}\|(L+n)u_n\|^2 <\infty.
\end{equation}

Along with $\bsymb{\CD}_\a$, we need its subset $\bsymb{\CD}_\a^0$
consisting of all elements $U\in \bsymb{\CD}_\a$ which have only a
finite number of non-zero components. Taking each component equal
zero in a vicinity of the point $x=0$, we obtain a dense subset in
$\GH$. Hence, both $\bsymb{\CD}_\a^0$ and $\bsymb{\CD}_\a$ are
dense in $\GH$.

Define the operators $\BA_\a$ and $\BA_\a^0$ as
\begin{gather*}
\BA_\a U=\CA U\sim\{(L+n) u_n\},\qquad   Dom\, \BA_\a=\bsymb{\CD}_\a;\\
\BA_\a^0= \BA_\a\res\bsymb{\CD}_\a^0.
\end{gather*}
Evidently, the operator $\BA_\a^0$ is symmetric in $\GH$. Our goal
is to prove the following result.
\begin{thm}\label{sa.t1}
The operator $\BA_\a$ is self-adjoint and coincides with the
closure of $\BA_\a^0$.
\end{thm}
Note that for $\a\neq\sqrt 2$ the result is covered by \cite{sol},
Theorem 5.1. Nevertheless, below we give the full proof of
Theorem. We do this mostly in order to illustrate the usage of
Proposition \ref{aux.bia}. In \cite{sol} another, more
sophisticated technical tools were used for the proof.
\begin{proof}
First of all, we show that
\begin{equation}\label{sa.2}
\BA_\a=\bigl(\BA_\a^0\bigr)^*.
\end{equation}

The inclusion ``$\subset$'' in \eqref{sa.2} can be easily checked
by the direct inspection. To prove the reverse inclusion, suppose
that $V\sim\{v_n\}\in Dom\,(\BA_\a^0\bigr)^*$ and
$(\BA_\a^0\bigr)^*V=W\sim\{w_n\}$. According to the definition of
the adjoint operator, this means that for any
$U\sim\{u_n\}\in\bsymb{\CD}_\a^0$ we have
\begin{equation}\label{sa.3}
\sum_{n\in\N_0}\int_{\R}(L+n)u_n\overline{v_n}dx=
 \sum_{n\in\N_0} \int_{\R} u_n\overline{w_n}dx.
\end{equation}
Take a function $f\in\CW(\R)$ such that $f(x)=0$ in a vicinity of
$x=0$ and fix a number $n_0\in\N_0$. The element $U\sim\{u_n\}$,
such that $u_{n_0}=f$ and $u_n=0$ for $n\neq n_0$, belongs to
$\bsymb{\CD}_\a^0$. Applying the identity \eqref{sa.3} to all such
$U$, we conclude that if $V\in Dom\, (\BA_\a^0\bigr)^*$, then
$v_n\in\CW(\R)$ and $w_n=(L+n)v_n$ for all $n\in\N_0$. Hence, for
$\{v_n\}$ the condition \eqref{sa.1} is satisfied. It remains to
check that the matching conditions \eqref{prel.16} are also
fulfilled. To this end, we fix a number $n_0\in\N_0$ and choose an
element $U\sim\{u_n\}$ as follows: $u_n(x)\equiv 0$ for
$|n-n_0|>1$; the functions $u_{n_0}, u_{n_0\pm 1}$ are supported
in a vicinity of the $x=0$ and in some smaller vicinity are given
by $u_{n_0}(x)=1$; $u_{n_0\pm 1}(x)=h_{\pm}|x|$ where $h_{\pm}$
are some appropriate numbers. Then $U\in\bsymb{\CD}_\a^0$,
provided that
\begin{equation*}
2\mu h_+=\sqrt{n_0+1},\qquad 2\mu h_-=\sqrt{n_0}.
\end{equation*}
Now \eqref{prel.16} for $n=n_0$ is implied by \eqref{sa.3} for the
element $U$ constructed. So, the equality \eqref{sa.2} is
justified. \vskip0.2cm

The statement of Theorem is equivalent to the fact that both
deficiency indexes of the operator $\BA_\a^0$ are equal to zero.
Since all the coefficients in the equation and in the matching
conditions are real, it is enough to prove that the only solution
$V\in\GH$ of the equation
\begin{equation}\label{sa.5}
\BA_\a V-iV=0
\end{equation}
is $V\equiv0$. Equation \eqref{sa.5} is nothing but \eqref{aux.1}
for $\L=i$. Using the representation \eqref{aux.equ}, we reduce
the equation to the form \eqref{aux.difeq} and can apply Lemma
\ref{sys.lem}.

If $\mu<1$, then according to \eqref{res.r3} the system  has a
pair of linearly independent solutions $\{C^\pm_n\}$ such that
\begin{gather}
|C^{\pm}_n|^2 \sim
n^{-1\pm\frac{\mu}{\sqrt{1-\mu^2}}}.\label{sa.di}
\end{gather}
In view of \eqref{A0.nnorm}, only the sequence $\{C_n^-\}$ may
generate a solution $V\in\GH$ of equation \eqref{sa.5}.

Equation \eqref{aux.difeq} is of the form \eqref{aux.rec}. Now we
use for it the identity \eqref{aux.z}, with $C_n=C^{-}_n$. This
gives
\begin{equation}\label{sa.14x}
\sum_{n=0}^N|C^-_n|^2\im
\z_n=-\sqrt{N+1}\im(C^-_{N+1}\overline{C^-_N}).
\end{equation}
By \eqref{sa.di}, the right-hand side of \eqref{sa.14x} vanishes
as $N\to\infty$. Since $\im \z_n>0$ for each $n$, we conclude from
\eqref{sa.14x} that $V\equiv0$. This shows that for $\a>\sqrt2$
the operator $\BA_\a$ is self-adjoint.

\vskip0.2cm If $\a<\sqrt 2$, then $\mu>1$ and according to
\eqref{res.r3} one of the solutions  of the system
\eqref{aux.difeq} exponentially grows and another exponentially
decays. Only the latter may give rise to the solution $V\in\GH$ of
equation \eqref{sa.5}. Again, using the identity \eqref{aux.z} we
conclude that this solution is identically zero.

\vskip0.2cm Finally, let $\a=\sqrt 2$, then $\mu=1$ and the
formula \eqref{res.r3} gives
\begin{equation*}
C^\pm_n\sim e^{\pm\sqrt{2n}(1+i)}n^{-1/4}.
\end{equation*}
Only the sequence $\{C^-_n\}$ lies in $\ell^2$ and again, we
conclude from \eqref{aux.z} that the deficiency indexes are equal
to $0$.
\end{proof}

\section{Spectrum of the operators $\BA_\a$}\label{main}
\subsection{Results.}\label{r}
The following theorem is the central result of the paper.
\begin{thm}\label{main.abs}
Let $\a>0$, $\mu=\sqrt 2/\a$ and let $\CJ_0(\mu)$ be the Jacobi
matrix (operator), defined in \eqref{res.j0}. Then the a.c.
spectrum of the operator $\BA_\a$ and its multiplicity are
described by the equalities \eqref{int.1} and \eqref{int.2}. In
particular, for $\a<\sqrt2$
\begin{gather}
\begin{cases}\s_{a.c.}(\BA_\a)=[1/2,\infty),\\
\gm_{a.c.}(E;\BA_\a)=\gm_{a.c.}(E;\BA_0),\qquad a.e.\ E\ge 1/2.
\end{cases}\label{main.mult}\end{gather}
Further,
\begin{gather}
\s_{a.c.}(\BA_{\sqrt 2})=[0,\infty);\qquad
\s_{a.c.}(\BA_{\a})=\R,\ \a>\sqrt 2,\label{main.big}
\end{gather}
and for $\a\ge\sqrt 2$ and $E\in\s_{a.c.}(\BA_{\a})$ we have
\begin{equation}\label{main.mubig}
\gm_{a.c.}(E;\BA_\a)=\gm_{a.c.}(E;\BA_0)+1.
\end{equation}
\end{thm}
The equalities \eqref{main.mult}, \eqref{main.big}
 and \eqref{main.mubig} immediately follow
from the relations \eqref{int.1}, \eqref{int.2} and
\thmref{res.j0spec}. So, our main goal for the rest of the paper
is to prove \eqref{int.1} and \eqref{int.2}.

\vskip0.2cm

To make the picture more complete, we present also the result
concerning the point spectrum $\s_p(\BA_\a)$. Within minor detail,
these results were proved in the papers \cite{sol}, \cite{sol1}.

\begin{thm}\label{main.p}
1. For any $\a\ge0$ the operator $\BA_\a$ has no eigenvalues
$E\ge1/2$, and $\s_p(\BA_\a)=\emptyset$ for $\a\ge\sqrt2$.

2. For $0<\a<\sqrt 2$ the operator $\BA_\a$ is positive definite.
Its point spectrum is always non-empty and finite, and the number
$N(1/2;\BA_\a)$ of its eigenvalues (counted with multiplicities)
satisfies the asymptotic formula
\begin{equation*}
N(1/2;\BA_\a)\sim\frac1{4\sqrt{2(\mu(\a)-1)}},\qquad
\a\nearrow\sqrt2\qquad (\mu(\a)=\frac{\sqrt2}{\sqrt\a}).
\end{equation*}
\end{thm}

It follows from Theorems \ref{main.abs} and \ref{main.p} that the
essential and the absolutely continuous spectra of the operator
$\BA_\a$ coincide as sets. However, our approach does not show
that $\BA_\a$ has no singular continuous spectrum.

\subsection{Outline of proof of \thmref{main.p}.}\label{prf}
1. If $E\ge1/2,\ V\sim\{v_n\}$ and $\BA_\a V=E V$, then each
component $v_n$ lies in $\CW(\R)$ and satisfies equation
\eqref{aux.w} with $\L=E$. For $n\le E-1/2$ this yields $v_n\equiv
0$. This can be interpreted as the equality $C_n=0$ for the
coefficients in \eqref{aux.difeq}. Then these equations imply that
$C_n=0$ also for all $n>E-1/2$. \vskip0.2cm

The absence of eigenvalues $E<1/2$ for the operator $\BA_\a$ with
$\a>\sqrt2$ was proved in \cite{sol}, Theorem 7.1. The possibility
to use Proposition \ref{aux.bia} simplifies the proof, and also
allows one to include the borderline case $\a\ge\sqrt2$. We leave
this to the reader.

The statement 2 is covered by \cite{sol}, Theorem 6.2 and
\cite{sol1}, formula (3.10).

\section{Representation of the resolvent}\label{res}
\subsection{Auxiliary considerations}\label{res.pre}
In this section we derive a convenient representation of the
operator
\begin{equation*}
(\BA_\a-\L)^{-1}-(\BA_0-\L)^{-1},\qquad \L\notin\R.
\end{equation*}
The equality \eqref{res.repr} which we establish in
\thmref{res.t1} can be interpreted in terms of the extension
theory of symmetric operators. However, formally we do not use
this theory in our construction. \vskip0.2cm

Given an element $F\sim\{f_n\}\in\GH$, denote
\begin{equation*}
U_\a\sim\{u_{\a,n}\}=(\BA_\a-\L)^{-1}F,\qquad \a\ge0.
\end{equation*}
The functions $u_{0,n}$ satisfy the equation
\begin{equation*}
  -u_{0,n}''+(n+1/2-\L)u_{0,n}=f_n
\end{equation*}
and lie in $\tH^2(\R)$. Solving this equation, we find that
\begin{equation*}
2\z_n(\L) u_{0,n}(x)=\int_\R e^{-\z_n(\L)|x-t|}f_n(t)dt.
\end{equation*}
Denote
\begin{equation}\label{res.Jn}
J_n=\int_\R\y_n(t;\L)
f_n(t)dt=\bigl(f_n,\y_n(.;\overline\L)\bigr),
\end{equation}
then
\begin{equation}\label{res.Jna}
2\z_n(\L)(n+1/2)^{1/4} u_{0,n}(0)=J_n,\qquad n\in\N_0.
\end{equation}
Recall that the numbers $\z_n(\L)$ and the functions $\y_n(.;\L)$
were defined in \eqref{A0.zn} and \eqref{A0.yn}.

Let now $U_\a-U_0\sim\{v_n\}$. Each function $v_n$ satisfies the
homogeneous equation \eqref{aux.w} and belongs to $\CW(\R)$ and
hence, $v_n(x)=C_n\y_n(x;\L)$. The coefficients $C_n$ are
determined by the matching conditions for $U_\a$.  Since in view
of \eqref{res.Jna}
\begin{gather}
(n+1/2)^{1/4}C_n=v_n(0)=u_{\a,n}(0)-u_{0,n}(0)\label{res.Xn}\\=u_{\a,n}(0)-
(2\z_n(\L)(n+1/2)^{1/4})^{-1}J_n\notag
\end{gather}
and the derivative $u'_{0,n}$ is continuous at $x=0$, we get
\begin{gather*}
u'_{\a,n}(0+)-u'_{\a,n}(0-)=v'_n(0+)-v'_n(0-)\\ =
-2\z_n(\L)(n+1/2)^{1/4}C_n=(n+1/2)^{-1/4}J_n-2\z_n(\L)
u_{\a,n}(0).
\end{gather*}

It is convenient for us to denote
\begin{equation*}
X_n=(n+1/2)^{-1/4}u_{\a,n}(0),
\end{equation*}
then
\begin{gather*}
u'_{\a,n}(0+)-u'_{\a,n}(0-)=(n+1/2)^{-1/4}J_n-2\z_n(\L)(n+1/2)^{1/4}X_n
\end{gather*}
and the matching conditions \eqref{prel.16} reduce to
\begin{gather*}
\mu\bigl((n+1/2)^{-1/4}J_n-2\z_n(\L)(n+1/2)^{1/4}X_n\bigr)\\=
(n+1)^{1/2}(n+3/2)^{1/4}X_{n+1}+n^{1/2}(n-1/2)^{1/4}X_{n-1},
\end{gather*}
or
\begin{gather*}
d_{n+1}X_{n+1}+2\mu\, y_n(\L)X_n+d_nX_{n-1}=
\mu\bigl(f_n,{\y_n}(.;\overline{\L})\bigr).
\end{gather*}
This is the non-homogeneous counterpart of the recurrence system
\eqref{aux.difeq}. It can be re-written in terms of the matrix
$\CJ(\L;\mu)$ introduced in \eqref{res.j}:
\begin{equation}\label{res.Y}
\CJ(\L;\mu)X=\mu\bigl(f_n,{\y_n}(.;\overline{\L})\bigr),\qquad
X=\{X_n\}.
\end{equation}

\subsection{Basic formula}\label{res.basic}
For $\L\notin\R$ consider the operator
\begin{equation}\label{res.bt}
\BT(\L):\ell^2(\N_0)\to\GH,\qquad\BT(\L)\{X_n\}\sim
\{X_n{\y_n}(.;\L)\},
\end{equation}
by \eqref{A0.nnorm} it is bounded and has bounded inverse. Its
adjoint acts from $\GH$ to $\ell^2$ as
\begin{equation*}
\BT(\L)^*F=\bigl\{\int_\R f_n(x){\y_n}(x;\overline\L)dx\},\qquad
F\sim\{f_n\}.
\end{equation*}
Substituting $\overline\L$ for $\L$, we obtain that
\begin{equation*}
\BT(\overline{\L})^*F=\bigl\{\bigl(f_n,{\y_n}(.;\L)\bigr)\bigr\}.
\end{equation*}
Both operator-valued functions $\BT(\L)$ and
$\BT(\overline{\L})^*$ are analytic in $\C_\pm$.

Now we are in a position to present the basic formula which
relates the operator $\BA_\a$ to the Jacobi matrix $\CJ(\L;\mu),\
\mu=\sqrt2/\a$.
\begin{thm}\label{res.t1}
Let $\a>0,\ \mu=\sqrt2/\a$ and $\L\notin\R$. Then
\begin{gather}
(\BA_\a-\L)^{-1}-(\BA_0-\L)^{-1} =\BT(\L)\bigl(\mu
\CJ(\L;\mu)^{-1}-(2\tY(\L))^{-1}\bigr)\BT(\overline{\L})^*.\label{res.repr}
\end{gather}
\end{thm}
\begin{proof}
By \eqref{res.Jn} and \eqref{res.Xn}, we have
$C_n=X_n-(2y_n(\L))^{-1}\bigl(f_n,{\y_n}(.;\overline\L)\bigr)$, or
\begin{equation*}
\{C_n\}=X-(2\tY(\L))^{-1}\BT(\overline\L)^*,
\end{equation*}
whence
\begin{equation*}
U_\a-U_0\sim\{v_n\}=\BT(\L)\{C_n\}.
\end{equation*}
We also find from \eqref{res.Y} that
$X=\mu\CJ(\L;\mu)^{-1}\BT(\overline\L)^*$. The desired equality
\eqref{res.repr} is an immediate consequence of the three last
equations.
\end{proof}

\section{Proof of \thmref{main.abs}: $E<1/2$}\label{ml}
\subsection{Absolutely continuous spectrum and jump of the bordered resolvent.}
\label{ml.ac} Let $\Om\subset\C$ be a domain, symmetric with
respect to the real axis and containing an interval $I\subset\R$.
Let $\BZ(\L)$ be an operator-valued function which is analytic in
$\Om\cap\C_+$ and in $\Om\cap\C_-$. Its jump at a point $E\in I$
is defined as
\begin{gather*}
\jp\BZ\pj(E):=\jp\BZ(\L)\pj_{\L\to E+i0}=(s){\text{-}}\lim_{\L\to
E+i0}\bigl(\BZ(\L)-\BZ(\bar{\L})\bigr)\\ ({\text{non-tangential
limit}}),
\end{gather*}
provided that the (strong) limit does exist.

One of the ways to investigate the absolutely continuous spectrum
of a self-adjoint operator $\BA$ consists in studying the jump of
its "bordered resolvent" $\BZ_\BG(\L):=\BG^*(\BA-\L)^{-1}\BG$.
Here $\BG$ is an appropriate bounded operator. It is well known
that for any $\BA$ and any $\BG\in\GS_2$ the non-tangential limits
of the bordered resolvent $\BZ_\BG(\L)$ as $\L\to E\pm i0$ exist
a.e. and belong to $\GS_1$, see e.g. \cite{BE}. Hence, the jump
also exists a.e. (with respect to the Lebesgue measure). For a
given operator $\BA$ the jump of its bordered resolvent may exist
also for a wider class of borderings.

The following statement of a rather general nature plays the key
role in our analysis. Its proof is given in Appendix B.
\begin{thm}\label{Aa.l1b}
Let $\BA$ be a self-adjoint operator in a separable Hilbert space
$\CH$ and $\D\subset\R$ a given Borelian subset. Let $\BG$ be a
bounded operator with the dense range, such that the operator
\begin{equation}\label{Aa.jum}
    \BV_{\BA,\BG}(E)=(s){\text{-}}\lim_{\L\to E+i0}\BG^*\bigl((\BA-\L)^{-1}-
    (\BA-\overline\L)^{-1}\bigr)\BG
\end{equation}
is well-defined a.e. on $\D$. Then
\begin{equation}\label{Aa.dim}
\gm_{a.c.}(E;\BA)=\rank \BV_{\BA,\BG}(E) \qquad {\text{a.e. on}}\
\D.
\end{equation}
\end{thm}
\begin{cor}\label{Aa.cor}
Let the assumptions of \thmref{Aa.l1b} be fulfilled for an
interval $\D\subseteq\R$. If $\rank \BV_{\BA,\BG}(E)>0$ a.e. on
$\D$, then $\D\subset \s_{a.c.}(\BA)$, and if $\rank
\BV_{\BA,\BG}(E)=0$ a.e. on $\D$, then
$\s_{a.c.}(\BA)\cap\D=\emptyset$.
\end{cor}

\subsection{Bordered resolvent of the operator
$\BA_\a$}\label{Aa.br} Along with the operator
$\tY_0=\CD\{n+1/2\}$ acting in $\ell^2$, let us define
$\GY_0=\BI\otimes\tY_0$. The latter operator acts in the space
$\GH$ interpreted as the tensor product $\ell^2\otimes\tL^2(\R)$,
cf. Section \ref{prel}. The powers $\tY_0^{-\g},\ \g>0$ are
compact operators in $\ell^2$, while the powers $\GY_0^{-\g}$ are
only bounded operators. If $U\sim\{u_n\}\in\GH$, then
$\GY_0^{-\g}U\sim\{(n+1/2)^{-\g}u_n\}$. Note that
\begin{equation}\label{Aa.js}
\GY^{-\g}\BT(\L)=\BT(\L)\tY^{-\g},\qquad
\BT(\overline\L)^*\GY^{-\g}= \tY^{-\g}\BT(\overline\L)^*.
\end{equation}
We shall show that the boundary limits of the bordered resolvent
$\BG^*(\BA_\a-\L)^{-1}\BG$ do exist if we take $\BG=\GY_0^{-\g}$
with $\g>1/4$, though $\BG$ is non-compact, and even the operator
$\tY_0^{-\g}$ lies in $\GS_2$ only if $\g>1/2$.

Let us consider the operator-valued functions
\begin{gather}
\BZ_\g(\L;\mu)=\GY_0^{-\g}(\BA_\a-\L)^{-1}\GY_0^{-\g},\qquad
\a=\sqrt2/\mu,
\notag\\
\wt{\BZ_\g}(\L;\mu)=\tY_0^{-\g}\CJ(\L;\mu)^{-1}\tY_0^{-\g},\notag\\
\wt{\BZ_{0,\g}}(\L;\mu)=\tY_0^{-\g}(\CJ_0(\mu)-\mu\L)^{-1}\tY_0^{-\g}.
\label{Aa.jJ0}
\end{gather}
In \eqref{Aa.jJ0} $\CJ_0(\mu)$ is the Jacobi matrix introduced in
\eqref{res.j0}.

By \eqref{res.repr} and \eqref{Aa.js},
\begin{gather*}
\BZ_\g(\L;\mu)=\GY_0^{-\g}(\BA_0-\L)^{-1}\GY_0^{-\g}\\+\BT(\L)\biggl(
\mu\wt{\BZ_\g}(\L;\mu)-\tY_0^{-\g}(2\tY(\L))^{-1}\tY_0^{-\g}\biggr)
\BT(\overline{\L})^*.
\end{gather*}
The jumps at a point $E\in\R$ of all but one terms in the
right-hand side are evidently equal to zero, provided that
$E<1/2$, and we get, at least formally,
\begin{equation}\label{Aa.jump}
\jp\BZ_\g(\L;\mu)\pj_{\L\to
E+i0}=\mu\BT(E)\jp\wt{\BZ_\g}(\L;\mu)\pj_{\L\to E+i0}\BT(E)^*,\
E<1/2.
\end{equation}
Taking into account the equalities \eqref{res.Psi} and
\eqref{res.psi}, it is natural to expect that the jumps
$\jp\wt{\BZ_\g}(.;\mu)\pj(E)$ and
$\jp\wt{\BZ_{0,\g}}(.;\mu)\pj(E)$ are close to each other.
Together with \eqref{Aa.jump}, this would allow us to express the
quantity $\jp\BZ_\g(.;\mu)\pj(E)$ through
$\jp\wt{\BZ_{0,\g}}(.;\mu)\pj(E)$, which makes it possible to use
\thmref{Aa.l1b}. Now we proceed to the successive realization of
this program.

\bigskip

We start with the following lemma.
\begin{lem}\label{Aa.l1}
Let $\mu> 0$, $\L\neq\bar\L\in\C$ and $\g>1/4$. Then
$\wt{\BZ_\g}(\L;\mu)\in\GS_1$. Besides, the strong non-tangential
limits $\wt{\BZ_\g}(E\pm i0;\mu)$ exist for almost all $E<1/2$.

The same results are valid for $\wt{\BZ_{0,\g}}(\L;\mu)$.
\end{lem}
\begin{proof} For definiteness, we consider $\L\in\C_+$.
The operator-valued function $\wt{\BZ_\g}(\L;\mu)$ is analytic in
$\L$ in the half-plane $\C_+$. Since $\im \CJ(\L;\mu)\le 0$, we
conclude that $\im \wt{\BZ_\g}(\L;\mu)\ge 0$, i.e. the values of
this function are dissipative operators. By a result due to
S.Naboko, see \cite{Na0}, Remark (1) to Theorem 2.2, any such
function admits the representation
\begin{equation}\label{Aa.rep}
\wt{\BZ_\g}(\L;\mu)=A+B\L+R^*(I+\L\CL)(\CL-\L)^{-1}R,
\end{equation}
where $\CL$ is a self-adjoint operator in an auxiliary Hilbert
space $\GH_0$, $A=A^*$ and $B=B^*\ge 0$ are bounded operators in
$\ell^2$, and $R$ is a bounded operator from $\ell^2$ to $\GH_0$.
It immediately follows from \eqref{res.norm} that in this
representation $B=0$. The operators $\CL,\, A$ and $R$ may depend
on $\mu$ but we do not reflect this dependence in our notations.

Taking in \eqref{Aa.rep} $\L=i$, we find that
\begin{equation}\label{Aa.re}
\wt{\BZ_\g}(i;\mu)=A+iR^*R.
\end{equation}
Further, let us show that for any $\g>0$ one has
$\wt{\BZ_\g}(i;\mu)\in\GS_2$. To this end, consider the inverse
matrix $\{G_{n,k}\}:=\CJ(i;\mu)^{-1}$. This matrix can be
expressed through the solutions of the homogeneous system
\eqref{aux.difeq} for $\L=i$. Namely, it is symmetric (i.e.
$G_{n,k}=G_{k,n}$) and its entries for $k\ge n$ are
\begin{equation*}
G_{n,k}=C g^{(1)}(k)g^{(2)}(n)
\end{equation*}
where $C$ is a constant, $\{g^{(1)}(k)\}$ is the decaying solution
of the system, and $\{g^{(2)}(n)\}$ is some other solution, a
certain linear combination of two basic solutions. The matrix
$\wt{\BZ_\g}(i;\mu)$ has the entries
\begin{equation*}
    (n+1/2)^{-\g}G_{n,k}(k+1/2)^{-\g},
\end{equation*}
and the standard calculation shows that
$\wt{\BZ_\g}(i;\mu)\in\GS_2$ for any $\g>0$ and an arbitrary
$\mu>0$. This immediately yields that
$\wt{\BZ_\g}(i;\mu)\in\GS_1$, provided that $\g>1/4$.

Now we conclude from \eqref{Aa.re} that $R\in\GS_2$. Hence,
$\wt{\BZ_\g}(\L;\mu)\in\GS_1$ for all $\L\in\C_+$. Besides, by a
result of \cite{BE}, Lemma 2.4, this implies the existence a.e. of
the non-tangential strong limits of the bordered resolvent
$R^*(\CL-\L)^{-1}R$. The equality \eqref{Aa.rep} shows that such
limits do exist also for the operator-valued function
$\wt{\BZ_\g}(\L;\mu)$.

For the operator-valued function $\wt{\BZ_{0,\g}}$ the proof is
simpler, since we need not the representation \eqref{Aa.rep}.
Otherwise, the argument remains the same.
\end{proof}

As a result of Lemma \ref{Aa.l1}, the equality \eqref{Aa.jump} is
justified.

\subsection{Jumps of $\wt{\BZ_\g}$ and of
$\wt{\BZ_{0,\g}}$.}\label{Aa.jumps} A subset $\D
\subset(-\infty,1/2)$ of full measure can be selected, such that
both these jumps are well-defined for $E\in\D$. According to the
Hilbert identity and using the equality \eqref{res.Psi}, we find
that
\begin{gather}
\wt{\BZ_\g}(\L;\mu)-\wt{\BZ_{0,\g}}(\L;\mu)=\tY_0^{-\g}\bigl(\CJ(\L;\mu)^{-1}-
(\CJ_0(\mu)-\mu\L)^{-1}\bigr)\tY_0^{-\g}\label{Aa.hil}\\
=-\tY_0^{-\g}\CJ(\L;\mu)^{-1}\Psi(\L;\mu)(\CJ_0(\mu)-\mu\L)^{-1}\tY_0^{-\g}
\notag\\=-\wt{\BZ_\g}(\L;\mu)\Phi(\L;\mu)\wt{\BZ_{0,\g}}(\L;\mu)\notag
\end{gather}
where
\begin{equation}\label{Aa.Phi}
\Phi(\L;\mu)=\tY_0^{\g}\Psi(\L;\mu)\tY_0^{\g}.
\end{equation}
More exactly, $\Phi(\L;\mu)$ is the extension by continuity of
this operator, defined originally on the set $Ran\, \tY_0^\g$.

The calculations below are carried through for a fixed value of
$\mu$, and we drop this parameter from our notations. It follows
from \eqref{res.psi} that for $\g<1/2$ the operator-valued
function $\Phi(\cdot)$ is analytic in $\C_\pm$ and its values are
compact operators in $\ell^2$. Below we always assume
$1/4<\g<1/2$, then both operator-valued functions
$\wt{\BZ_\g}(\L)$ and $\wt{\BZ_{0,\g}}(\L)$ have boundary limits
as $\L\to E\pm i0$ a.e. on $(-\infty,1/2)$.

The equality \eqref{Aa.hil} and the similar equality for $\bar\L$
imply that
\begin{gather*}
\wt{\BZ_{0,\g}}(\L)-\wt{\BZ_{0,\g}}(\bar\L)=\wt{\BZ_\g}(\L)-
\wt{\BZ_\g}(\bar\L)\\
+\bigl(\wt{\BZ_\g}(\L)-\wt{\BZ_\g}(\bar\L)\bigr)\Phi(\L)
\wt{\BZ_{0,\g}}(\L)
+\wt{\BZ_\g}(\bar\L)\bigl(\Phi(\L)-\Phi(\bar\L)\bigr)\wt{\BZ_{0,\g}}(\L)\\
+\wt{\BZ_\g}(\bar\L)\Phi(\bar\L)\bigl(\wt{\BZ_{0,\g}}(\L)-
\wt{\BZ_{0,\g}}(\bar\L)\bigr).
\end{gather*}
Letting $\L\to E+i0$ non-tangentially, we obtain for a.e. $E<1/2$:
\begin{gather}
\jp\wt{\BZ_{0,\g}}\pj(E)=\jp\wt{\BZ_{\g}}\pj(E)+
\jp\wt{\BZ_{\g}}\pj(E)\Phi(E)
\wt{\BZ_{0,\g}}(E+i0)\label{Aa.jps}\\
+ \wt{\BZ_\g}(E-i0)\Phi(E)\jp\wt{\BZ_{0,\g}}\pj(E).\notag
\end{gather}

Let us consider the operator-valued functions
\begin{gather}
\CG_+(\L)=I+\Phi(\L)\wt{\BZ_{0,\g}}(\L),\qquad
\L\in\C_+,\label{Aa.gpl}\\
\CG_-(\L)=I-\wt{\BZ_\g}(\L)\Phi(\L),\qquad
\L\in\C_-.\label{Aa.gmin}
\end{gather}
The function $\CG_+(\L)$ can be represented in a different way:
using the definitions \eqref{Aa.Phi} and \eqref{Aa.jJ0}, we find
that
\begin{gather}
\CG_+(\L)=I+\tY_0^{\g}\bigl(\CJ(\L)-(\CJ_0-\mu\L)\bigr)\tY_0^{\g}\tY_0^{-\g}
(\CJ_0-\mu\L)^{-1}\tY_0^{-\g}\notag\\
=\tY_0^{\g}\CJ(\L)(\CJ_0-\mu\L)^{-1}\tY_0^{-\g}.\label{Aa.gpl1}
\end{gather}
This calculation shows also that the right-hand side in
\eqref{Aa.gpl1} is well-defined as a bounded operator in $\ell^2$.
For the "minus" sign, it is more convenient to deal with the
adjoint operator, and we get
\begin{equation}\label{Aa.gmin1}
\CG_-(\L)^*=\tY_0^{\g}(\CJ_0-\mu\overline\L)\CJ(\overline\L)^{-1}\tY_0^{-\g}.
\end{equation}

It follows from the definitions \eqref{Aa.gpl} and \eqref{Aa.gmin}
and from compactness of $\Phi(\L)$ that the functions
$\CG_\pm(\L)$ are analytic in $\C_\pm$ respectively, and for each
$\L\in\C_\pm$ the operator $\CG_\pm(\L)-I$ is compact. Hence, the
image of $\CG_\pm(\L)$ is a closed subset in $\ell^2$.

Now, we need the following Lemma.
\begin{lem}\label{Aa.l2}
For any $\L\in\C_\pm$ the operator $\CG_\pm(\L)$ has bounded
inverse.
\end{lem}
\begin{proof}
Throughout the proof, $(.,.)$ stands for the scalar product in
$\ell^2$. By \eqref{Aa.gpl1}, we have for any $f\in\ell^2$ and any
$g\in Dom\, \tY_0^\g$:
\begin{equation}\label{Aa.aux}
\bigl(\CG_+(\L)f,g\bigr)=\bigl(\CJ(\L)(\CJ_0-\mu\L)^{-1}
\tY_0^{-\g}f,\tY_0^{\g}g\bigr).
\end{equation}
Suppose that $f\in\nul\,\CG_+(\L)$ for some $f\neq 0$ and
$\L\in\C_+$. Since $Ran\, \tY_0^\g=\ell^2$, we conclude from
\eqref{Aa.aux} that then $\nul\, \CJ(\L)$ is non-trivial which
contradicts the dissipativity of $\CJ(\L)$.

In a similar way, we derive from \eqref{Aa.gmin1} that for $f\in
Dom\,\tY_0^\g$ and any $g\in\ell^2$
\begin{equation*}
\bigl(\CG_-(\L)f,g\bigr)=\bigl(\tY_0^\g
f,(\CJ_0-\mu\overline\L)\CJ(\overline\L)^{-1}\tY_0^{-\g}g\bigr).
\end{equation*}
Suppose that $Ran \CG_-(\L)\neq\ell^2$ for some $\L\in\C_-$. Then
there exists an element $g\in\ell^2$, such that
$(\CG_-(\L)f,g\bigr)=0$ for all $f\in\ell^2$. This would imply
\begin{equation*}
(\CJ_0-\mu\bar\L)\CJ(\bar\L)^{-1}\tY_0^{-\g}g=0.
\end{equation*}
However, for $g\neq 0$ this is impossible, since
$\CJ(\bar\L)^{-1}$ and $\tY_0^{-\g}$ are inverse operators and
hence, have trivial kernels. Further,
$\nul(\CJ_0-\mu\bar\L)=\{0\}$ for $\L\notin\R$, since the operator
$\CJ_0$ is self-adjoint.
\end{proof}

Now we use the following statement, see \cite{Na}.
\begin{prop}\label{Aa.nab}
Let $\CX\in\C_\pm$ be a domain, such that $\overline{\CX}\cap\R$
contains an interval $\D$. Let $\CG(\L)$ be an analytic
operator-valued function in $\CX$, such that $\CG(\L)-I$ are
compact operators in a Hilbert space $\CH$. Suppose that for
almost all $E\in \D$ the function $\CG$ is non-tangentially
bounded at $E$ and has a strong non-tangential limit $\CG(E+i0)$,
such that $\CG(E+i0)-I\in \GS_\infty$. Suppose also that for at
least one point $\L_0\in\C_+$ the operator $\CG(\L_0)$ has bounded
inverse.

Then for almost all $E\in \D$ the operator $\CG(E+i0)$ has bounded
inverse.
\end{prop}

It follows from Lemmas \ref{Aa.l1} and \ref{Aa.l2} that the
assumptions of Proposition are fulfilled for the operator-valued
functions $\CG_\mp(\L)$ defined in \eqref{Aa.gpl} and
\eqref{Aa.gmin}. This allows us to conclude from \eqref{Aa.jump}
and \eqref{Aa.jps} that for a.e. $E<1/2$
\begin{equation}\label{Aa.5}
\jp\BZ_\g\pj(E)=\BT(E)\CG_-(E-i0)\jp\wt{\BZ_{0,\g}}\pj(E)\bigl(\CG_+
(E+i0)\bigr)^{-1}\BT(E)^*.
\end{equation}
The equalities \eqref{Aa.jump} and \eqref{Aa.5} yield that
\begin{equation*}
    \rank\jp\BZ_\g\pj(E)=\rank\jp\wt{\BZ_{0,\g}}\pj(E),\qquad
    a.e. \ E<1/2.
\end{equation*}
According to \thmref{Aa.l1b}, it follows that
\begin{equation*}
\gm_{a.c.}(E;\BA_\a)=\gm_{a.c.}(E;\CJ_0(\mu))\qquad a.e.\
{\text{on}}\ (-\infty,1/2).
\end{equation*}
Since $\gm_{a.c.}(E;\BA_0)=0$ for $E<1/2$, the equality
\eqref{int.2} for such $E$ is justified.

Now, we conclude from Corollary \ref{Aa.cor} that for $\a<\sqrt2$
the operator $\BA_\a$ has no a.c. spectrum below the point $1/2$.
For $\a>\sqrt2$ we have $(-\infty,1/2)\subset\s_{a.c.}(\BA_\a)$,
and $[0,1/2)\subset\s_{a.c.}(\BA_{\sqrt2})$.

\section{Proof of \thmref{main.abs}: $E\ge 1/2$}\label{big}
The proof extends to $\l>1/2$ with the help of a simple technical
trick. It is based on the passage to the subspace
$\GH_m=\ell^2(\N_m,\tL^2(\R))$, with $m$ large enough. Here
$\N_m=\{m,m+1,\ldots\}$, so that, in particular, $\N_1=\N$. The
subspace $\GH_m$ is not invariant for the operator $\BA_\a$,
however it is invariant for an appropriate operator
$\wh\BA^{(m)}_\a$, that can be obtained from $\BA_\a$ by the
perturbation of its resolvent by a finite rank operator. For the
operator $\wh\BA^{(m)}_\a$ the scheme developed in Section
\ref{ml} works for $E<m+1/2$, and return to the original operator
$\BA_\a$ does not change the absolutely continuous spectrum and
its multiplicity function.

\subsection{Operators $\BA^{(m)}_\a$ and $\wh\BA^{(m)}_\a$.}\label{Aan}
Denote $\N_m=\{m,m+1,\ldots\}$, so that, in particular, $\N_1=\N$.
Let us define an operator $\BA^{(m)}_\a$, acting in the space
$\GH_m=\ell^2(\N_m,\tL^2(\R))$. Namely, its domain
$\bsymb{\CD}^{(m)}_\a$ consists of the elements
$U\sim\{u_n\}_{n\ge m}$, such that each component $u_n$ lies in
$\CW(\R)$ (see \eqref{aux.CW}), the matching conditions
\begin{equation*}
u'_n(0+)-u'_n(0-)=\begin{cases}
\frac{\a}{\sqrt2}\sqrt{m+1}u_{m+1}(0),\qquad & n=m,\\
\frac{\a}{\sqrt2}
\bigl(\sqrt{n+1}u_{n+1}(0)+\sqrt{n}u_{n-1}(0)\bigr), \qquad & n>m
\end{cases}\end{equation*}
are satisfied, and $\sum_{n\ge m}\|(L+n)u_n\|^2 <\infty$. For
$U\in\bsymb{\CD}^{(m)}_\a$ we let
\begin{equation*}
 \BA^{(m)}_\a U\sim\{(L+n)u_n\}_{n\ge m},
\end{equation*}
cf. \eqref{prel.10}. Evidently, $\BA^{(0)}_\a=\BA_\a$. Each
operator $\BA^{(m)}_\a$ is self-adjoint, the proof is the same as
for $m=0$.

\vskip0.2cm Along with $\BA^{(m)}_\a$, we need also the operators
$\wh\BA^{(m)}_\a$, acting in the original Hilbert space $\GH$.
Namely, we let
\begin{equation}\label{res.3}
\wh\BA^{(m)}_\a=\biggl({\sum_{n<m}}^\oplus(\BH+n)\biggr)\oplus\BA^{(m)}_\a
\end{equation}
where $\BH$ is the operator in $\tL^2(\R)$, defined in
\eqref{A0.A}. Note that $\wh\BA^{(m)}_0=\BA_0$ for any $m$.
\vskip0.2cm Consider also two  sequences of Jacobi matrices,
$\CJ_0^{(m)}(\mu)$ and $\CJ^{(m)}(\L;\mu)$. Each
$\CJ_0^{(m)}(\mu)$ is a sub-matrix of the matrix $\CJ_0(\mu)$
given by \eqref{res.j0}, it is obtained from $\CJ_0(\mu)$ by
removing its first $m$ rows and $m$ columns. The matrix
$\CJ^{(m)}(\L;\mu)$ is obtained in the same way from the matrix
\eqref{res.j}.

We also define a sequence $\tY^{(m)}(\L)$ of diagonal operators in
$\ell^2(\N_m)$, namely $\tY^{(m)}(\L)=\CD\{y_n(\L)\}_{n\ge m}$. It
is clear that
\begin{equation*}
   \CJ_0^{(0)}(\mu)=\CJ_0(\mu), \ \CJ^{(0)}(\L;\mu)=\CJ(\L;\mu), \
\tY^{(0)}(\L)=\tY(\L).
\end{equation*}
Let $\BT^{(m)}(\L):\ell^2(\N_m)\to\GH_m$ be the natural
restriction of the operator $\BT(\L)$ defined in \eqref{res.bt}.
The equality
\begin{gather}
(\BA^{(m)}_\a-\L)^{-1}-(\BA^{(m)}_0-\L)^{-1}\label{res.4}\\
=\BT^{(m)}(\L)\bigl(\mu
\CJ^{(m)}(\L;\mu)^{-1}-(2\tY^{(m)}(\L))^{-1}\bigr)\BT^{(m)}
(\overline{\L})^*\notag
\end{gather}
can be justified in the same way as \eqref{res.repr}. It is
important that all the operator-valued functions appearing in
\eqref{res.4} are analytic in the domain
$\Om_m=\C\setminus[m+1/2,\infty)$.

\vskip0.2cm

It follows from \eqref{res.4} that the operator
\begin{equation*}
    (\wh\BA^{(m)}_\a-\L)^{-1}-(\BA_0-\L)^{-1}
\end{equation*}
is equal to the orthogonal sum of the null operator acting in the
space $\ell^2(\{0,\ldots,m-1\},\tL^2(\R))$ and the operator
\begin{equation*}
\BT^{(m)}(\L)\bigl(\mu\CJ^{(m)}(\L;\mu)^{-1}-(2\tY^{(m)}(\L))^{-1}\bigr)
\BT^{(m)}(\overline{\L})^*
\end{equation*}
The next statement is an immediate consequence of the above
reasonings.
\begin{lem}\label{res.le}
For any $\a>0$ and $m\in\N$ and for any $\L\notin[1/2,\infty)$ we
have
\begin{equation*}
   \rank\bigl( (\BA_\a-\L)^{-1}-(\wh\BA^{(m)}_\a-\L)^{-1}\bigr)=m
\end{equation*}
and therefore,
\begin{equation}\label{res.spec}
\s_{a.c.}(\wh\BA^{(m)}_\a)=\s_{a.c.}(\BA_\a),\qquad
\gm_{a.c.}(E;\wh\BA^{(m)}_\a)=\gm_{a.c.}(E;\BA_\a)\ a.e.
\end{equation}
\end{lem}

Now it is easy to conclude the proof of \thmref{main.abs} for
$E\ge 1/2$. Take $m\in\N$ such that $E<m+1/2$. \thmref{res.j0spec}
evidently applies to the matrices $\CJ_0^{(m)}$, and the scheme
developed in Section \ref{ml} works for the operator
$\wh\BA^{(m)}_\a$ without any change. One only has to keep in mind
that now in the corresponding version of Lemma \ref{Aa.l1} one can
take $E<m+1/2$. As a result, we obtain that
\begin{gather*}
\s_{a.c.}(\BA^{(m)}_\a)=\s_{a.c.}(\BA^{(m)}_0)\cup\s_{a.c.}(\CJ^{(m)}_0(\mu)),
\\
\gm_{a.c.}(E;\BA^{(m)}_\a)=\gm_{a.c.}(E;\BA^{(m)}_0)+\gm_{a.c.}(E;\CJ^{(m)}_0
(\mu))
\end{gather*}
Taking into account the equalities \eqref{res.3} and
\eqref{res.spec}, we arrive at the desired result.

\vskip0.2cm Note that for $\a<\sqrt2$ it is easy to prove
\thmref{main.abs} by means of the quadratic form approach, cf.
proof of Theorem 6.2 in \cite{sol}. However, in the present paper
we decided to give a unified exposition for all values of the
parameter.

\section{Concluding remarks}\label{conc}
In the model suggested by Smilansky in \cite{sm} the operators
$\BA_\a$ act in the space $\tL^2(\G\times\R)$ where $\G$ is a
metric graph (in another terminology, a quantum graph). The model
is interpreted as ``harmonic oscillator, attached to a graph''. In
order to describe the setting of this, more general problem,
consider first the case when $\G=\G_m$ is a star graph with $m$
bonds, each of infinite length. More precisely, $\G_m$ is the
union of $m$ half-lines $B_1,\ldots, B_m$, emanating from the
common vertex $o$, {\it the root} of the tree. Let
$t\in[0,\infty)$ stand for the coordinate along each bond. The
value $t=0$ corresponds to the root $o$. Each function $f$ on
$\G_m$ can be viewed as a family of $m$ functions
$f_j=f\bigl|_{B_j}$ defined on $[0,\infty)$. If each $f_j$ has the
derivative at $t=0$, we set
\begin{equation*}
[f'](o)=\sum_{j=1}^m f'_j(0).
\end{equation*}

In the case considered, the operator $\BA_\a$ in
$\tL^2(\G_m\times\R)$ is defined by the differential expression
\eqref{prel.eq} for $(x,q)\in B_j\otimes\R, \ j=1,\ldots,m$ and
the matching condition
\begin{equation}\label{conc.1}
[U'_x(.,q)](o)=\a qU(o,q),\qquad q\in\R.
\end{equation}
The real axis $\R$ with the marked point $o=0$ can be identified
with the graph $\G_2$, and in this case the condition
\eqref{conc.1} turns into \eqref{prel.tran}.

All the results of the present paper extend to the star graphs
$\G_m$ with an arbitrary $m>0$, with only minor changes: 1) the
equality in \eqref{A0.mulsac} has to be replaced by
$\gm_{a.c.}(E;\BA_0)=mn$; 2) the borderline point between the
small and the large values of $\a$ is now $m/\sqrt2$, and the
expression for $\mu$ becomes $\mu=m(\a\sqrt2)^{-1}$. The
equalities \eqref{int.1} and \eqref{int.2} survive. The proofs
basically remain the same, but the technical calculations
sometimes become rather lengthy. This was the only reason, why we
restricted ourselves to the case $m=2$ in the main part of this
paper.

\vskip0.2cm

In a similar way, the case when $\G$ is a general star graph with
$m$ bonds can be considered. Some of the bonds (say, $m_0$ where $
0\le m_0\le m$) are supposed to be of infinite, and other of
finite length. The Dirichlet boundary condition is imposed at the
ends of the finite bonds. The point spectrum and the essential
spectrum of the operators $\BA_\a$ were considered in \cite{sol}
for this case. The a.c. spectrum can be analyzed by means of the
same approach as in the present paper, but somewhat more serious
changes in the formulations are necessary. They stem from the fact
that now in the analog of \eqref{A0.mulsac} we have
$\gm_{a.c.}(E;\BA_0)=m_0n$. In particular, if $m_0=0$, that is if
the graph is compact, then the spectrum of $\BA_0$ is discrete. It
remains discrete for $\BA_\a$ with $\a<m/\sqrt2$, but its
absolutely continuous component fills the whole of $\R$ for $\a\ge
m/\sqrt2$.

\vskip0.2cm

The case of an arbitrary metric graph with a finite number of
bonds can be also analyzed, but this requires a bit more advanced
technical tools. Still, the main ideas remain the same. This
material will be presented elsewhere.

\appendix
\section{Proof of \thmref{res.j0spec}}\label{ap1}
The proof is based upon the Gilbert -- Pearson theory \cite{gp} of
subordinate solutions, or more exactly upon the version of this
theory for Jacobi matrices, see \cite{kp}.

First of all, we have to consider the homogeneous equation
\begin{equation*}
\CJ_0(\mu)h=z h,
\end{equation*}
or
\begin{equation}\label{a1}
d_{n+1}C(n+1)+\bigl((2n+1)\mu-z\bigr)C(n)+d_n C(n-1)=0.
\end{equation}
This equation is similar to \eqref{aux.difeq}, it can be written
in the form \eqref{aux.dieq}, and moreover, for the coefficients
in the decompositions \eqref{aux.dec} we have
\begin{equation*}
a_0=2\mu,\ a_1=-(\mu+z);\qquad b_0=1, \ b_1=-1.
\end{equation*}
Comparing this with \eqref{sys.1}, we conclude that the asymptotic
formulas \ref{res.r3} apply to the system  \eqref{a1} if we take
$\L=z/\mu$.

\vskip0.2cm

1. If $\mu>1$, the formula \eqref{res.r3} shows that for any
$E\in\R$ equation \eqref{a1} has a subordinate solution. Namely,
this is the solution with the sign ``$+$'' in the exponent. It
follows that the spectrum of $\CJ_0(\mu)$ is discrete. It is easy
to show that for $\mu>1$ the matrix $\CJ_0(\mu)$ is positive
definite, hence its eigenvalues tend to $+\infty$. \vskip0.2cm

2. If $\mu<1$, the formula \eqref{res.r3} shows that for any
$E\in\R$ equation \eqref{a1} has no subordinate solution, since
$|C^{\pm}(n)|\sim n^{-1/2}$.

As for any self-adjoint Jacobi matrix, the spectrum of
$\CJ_0(\mu)$ is simple, and the element $e_0=\{1,0,0\ldots\}^\top$
can be taken as the generating vector. Let $\CE$ stand for the
spectral measure of the operator $\CJ_0(\mu)$, then there exists a
non-negative function $\tau\in\tL^1(\R)$ such that
\begin{equation*}
\bigl(\CE(\D)e_0,e_0\bigr)_{\ell^2}=\int_\D\tau(\l)d\l
\end{equation*}
for an arbitrary Borelian set $\D\subset\R$. Our aim is to show
that $\tau(\l)\neq 0$ a.e. on $\R$.

Let $m_\infty(z)$ be the Weyl function for the equation
$\CJ(\mu)h-zh=0$. It is a Herglotz function, and therefore it has
boundary limits $m_\infty(E+i0)$ for a.e. $E\in\R$. It follows
from \cite{kp}, Theorem 1 that $\im m(E+i0)\neq 0$ a.e. According
to the formula (5) from \cite{kp}
$$
(\CJ_0(\mu)-z)^{-1}e_0,e_0)=\bigl(\mu-z+d_1m_\infty(z)\bigr)^{-1}.
$$
Recall that by \eqref{res.entr} $d_1=(3/4)^{1/4}$.

By the Spectral Theorem,
$$
(\CJ_0(\mu)-z)^{-1}e_0,e_0)=\int_\R\frac{\tau(\l)d\l}{\l-z}.
$$
The last two equalities imply (when $z=E+i\vare$ and $\vare\to+0$)
that
$$
\im\bigl(\mu-z+d_1m_\infty(z)\bigr)^{-1}\to \frac{-d_1\im
m_\infty(E+i0)} {|\mu-E+d_1m_\infty(E+i0)|^2}\neq0.
$$
This limit is equal to $\pi\tau(E)$. Therefore, $\tau(E)>0$ a.e.
which shows that $\gm(E;\CJ_0(\mu))=1$ a.e.

\vskip0.2cm

3. The matrix $\CJ_0(1)$ is non-negative, therefore its spectrum
lies on $[0,\infty)$. The formula \eqref{res.r3} shows that for
$E>0$ the equation $\CJ_0(1)h=Eh$ has no subordinate solution. The
equality $\gm_{a.c.}(E;\CJ_0(1))=1$ a.e. on $\R_+$ can be proved
in the same way as in the previous case.

This concludes the proof of \thmref{res.j0spec}.

\section{Proof of \thmref{Aa.l1b}}\label{ap2}
\subsection{Remarks on the a.c. spectrum and on its multiplicity.}\label{1}
Before giving the proof, we present some remarks of a rather
general nature, concerning the notion of a.c. spectrum. We
consider this useful because of the dual nature of the a.c.
spectrum. Indeed, it combines some features coming from the
measure theory with another ones, coming from topology of the real
line.

Let $\BK$ be a self-adjoint operator whose spectrum is purely
a.c., of multiplicity one. This means that $\BK$ is unitary
equivalent to the operator of multiplication, $u(\l)\mapsto\l
u(\l)$, in the space $\tL^2(\CX;d\l)$ where $\CX\subset\R$ is some
Borelian set and $d\l$ is the Lebesgue measure. The spectrum of
$\BK$ is the closure $\overline\CX$ of the set $\CX$, and it may
happen that $\meas(\overline\CX\setminus\CX)>0$. For this reason,
characterization of the a.c. spectrum of $\BK$ has to include
description of both the set $\s_{a.c.}(\BK)$ and the multiplicity
function $\gm_{a.c.}(\l;\BK)$. Say, in the above example
$\s_{a.c.}(\BK)=\overline\CX$ and the multiplicity function is
equal to one a.e. on $\CX$ and to zero a.e. outside $\CX$. It is
clear, how this extends to the case of spectrum of higher, or of
varying multiplicity.

Note that any Borelian set $\CX\subset\R$, such that
$\gm_{a.c.}(\l;\BK)>0$ a.e. on $\CX$ and  $\gm_{a.c.}(\l;\BK)=0$
a.e. on $\R\setminus\CX$, in the book \cite{Y} is called {\it the
core} of $\s_{a.c.}(\BK)$. \vskip0.2cm
\subsection{Preparatory material.}\label{2}
 We need some facts from the general theory of spectral
measure. Let $\BA$ be a self-adjoint operator in a separable
Hilbert space $\CH$ and $\CE$ be its spectral measure. For
elements $\phi,\psi\in\CH$, we denote by $\rho_{\phi,\psi}$ the
scalar complex-valued measure
$\rho_{\phi,\psi}(\cdot)=(\CE(\cdot)\phi,\psi)$ and by
$\rho'_{\phi,\psi}$ its Radon - Nikodym derivative with respect to
the Lebesgue measure. This derivative is defined a.e. on $\R$. It
is equal to zero if either of the elements $\phi,\psi$ is
orthogonal to the absolutely continuous subspace of the operator
$\BA$.

For any $\phi,\psi\in\CH$ there exists a subset
$\CX_{\phi,\psi}\subset\R$ of zero Lebesgue measure, such that
 \begin{equation}\label{ac.1}
\lim_{\L\to
E+i0}((\BA-\L)^{-1}\phi-(\BA-\overline\L)^{-1}\phi,\psi)= 2\pi
i\,\rho'_{\phi,\psi}(E),\qquad E\notin\CX_{\phi,\psi}.
\end{equation}
In such cases we do not use the notation $\jp\cdot\pj(E)$, since
the exceptional set $\CX_{\phi,\psi}$ depends on the chosen
elements. The set $\CX_{\phi,\psi}$ can be chosen in such a way
that the relation \eqref{ac.1} is satisfied for the pairs
$\{\phi,\phi\}$, $\{\psi,\psi\}$ and $\{\phi,\psi\}$
simultaneously, then it is satisfied also for any pair from the
linear hull of the elements $\phi$ and $\psi$.

\vskip0.2cm

Suppose now that $\BG$ is a bounded linear operator, such that the
jump \eqref{Aa.jum} exists a.e. on $\D$ where $\D\subset\R$ is a
given Borelian set. This implies that for any pair
$\phi,\psi\in\CH$ the limit
\begin{equation*}
\lim_{\L\to
E+i0}((\BA-\L)^{-1}\BG\phi-(\BA-\overline\L)^{-1}\BG\phi,\BG\psi)=
(\BV_{\BA,\BG}(E)\phi,\psi)
\end{equation*}
does exist for $E\notin\CX_\BG$ where $\CX_\BG\subset\D$ is a set
of Lebesgue measure $0$. Unlike \eqref{ac.1}, the set $\CX_\BG$
does not depend on the choice of $\phi$ and $\psi$. It follows
from \eqref{ac.1} that necessarily
\begin{equation}\label{Aa.1b}
(\BV_{\BA,\BG}(E)\phi,\psi)=2\pi i\,
\rho'_{\BG\phi,\BG\psi}(E)\qquad a.e.\ {\text{on}}\ \D.
\end{equation}
Again, here the exceptional set may depend on $\phi$ and $\psi$.
However, it is important that according to our assumption, for
$E\notin\CX_\BG$ the expression in the left-hand side of
\eqref{Aa.1b}, and thus the one in the right-hand side, is the
sesqui-linear form of a bounded operator.

\vskip0.2cm

Take a dense countable set $\{\phi_n\}$ of elements in $\CH$, then
there exists a subset $\CX'\subset\D$ of the Lebesgue measure $0$,
such that for $E\in\D\setminus\CX'$ the equality \eqref{Aa.1b} is
satisfied for all pairs $\{\phi_n,\phi_m\}$, and therefore for any
$\phi,\psi$ from the linear hull $\CM$ of the system
 $\{\phi_n\}$. In other words, {\it there exists a dense linear subspace
$\CM\subset\CH$, such that \eqref{Aa.1b} is satisfied for all
$E\in\D\setminus\CX'$ and for all $\phi,\psi\in\CM$
simultaneously}.

\subsection{Proof of the theorem.}\label{3}
Taking an appropriate partition of the original set, we may assume
that $\gm_{a.c.}(E;\BA)=\nu=const$ a.e. on $\D$. Suppose first
that $\nu<\infty$. According to the general theory of spectral
measure, see e.g. \cite{bs}, there is a subspace $\CH_\D
\subset\CH$ invariant with respect to $\BA$, isometric to
$\tL^2(\D,\ell^2_\nu;dx)$ and such that on $\CH_\D$ the operator
$\BA$ acts as multiplication by $x$. More precisely, let
$\Pi:\CH_\D\to \tL^2(\D,\ell^2_\nu;dx)$ stand for the above
isometry. Extending it by zero to the orthogonal complement of
$\CH_\D$, we obtain a partially isometric operator
$\Pi^\circ:\CH\to \tL^2(\D,\ell^2_\nu;dx)$. If $\phi\in
Dom\,\BA\cap\CH_\D$, then $\BA\phi\in\CH_\D$ and
\begin{equation*}
(\Pi\BA\phi)(x)=x(\Pi\phi)(x).
\end{equation*}

For any $\phi,\psi\in\CH_\D$ and any Borelian subset $\p\subset\D$
we have
\begin{equation*}
\rho_{\phi,\psi}(\p)=
\int_\p\bigl((\Pi\phi)(x),(\Pi\psi)(x)\bigr)_{\ell^2_\nu}dx
\end{equation*}
and
\begin{equation*}
\rho'_{\phi,\psi}(E)=\bigl((\Pi\phi)(E),(\Pi\psi)(E)\bigr)_{\ell^2_\nu},\qquad
a.e.\ {\text{on}}\ \D.
\end{equation*}
Moreover, $\rho'_{\phi,\psi}(E)=0$ a.e. on $\D$, provided that
either of the elements $\phi,\psi$ is orthogonal to $\CH_\D$. Now,
the formula \eqref{Aa.1b} and the above remarks imply that there
exists a dense linear subset $\CM\subset\CH$, such that for all
$\phi,\psi\in\CM$ the equality
\begin{equation}\label{Aa.bg}
(\BV_{\BA,\BG}(E)\phi,\psi)=2\pi i\,
\bigl((\Pi^\circ\BG\phi)(E),(\Pi^\circ\BG\psi)(E)\bigr)_{\ell^2_\nu}
\end{equation}
is satisfied a.e. on $\D$, and the exceptional subset does not
depend on the choice of $\phi,\psi$. The expression in the
right-hand side is the sesqui-linear form of an operator in
$\ell^2_\nu$. Necessarily, its rank does not exceed $\nu$.
Therefore, the restriction of the operator $\BV_{\BA,\BG}(E)$ to
$\CM$ also has rank no greater than $\nu$. Since the linear set
$\CM$ is dense in $\CH$ and the operator $\BV_{\BA,\BG}(E)$ is
bounded, we find that $\rank\BV_{\BA,\BG}(E)\le \nu$.

\vskip0.2cm

In order to obtain the reverse inequality, take elements
$f_1,\ldots,f_\nu\in\CH_\D$, such that
\begin{equation*}
    (\Pi f_j)(x)=\chi_\D(x)e_j,\qquad j=1,\ldots,\nu
\end{equation*}
where $\chi_\D$ is the characteristic function of the set $\D$ and
the vectors $\{e_j\}$ form the natural basis in $\ell^2_\nu$.
Since the range of $\BG$ is assumed dense, we can choose elements
$\phi_j\in\CH$ in such a way that
\begin{equation*}
\|\BG\phi_j-f_j\|_\CH<\vare, \qquad j=1,\ldots,\nu,
\end{equation*}
where $\vare>0$ is arbitrarily small. Then
\begin{equation*}
\int_\D\|(\Pi^\circ \BG\phi_j)(x)-\chi_\D(x)e_j\|^2_{\ell^2_\nu}dx
=\|\Pi^\circ(\BG\phi_j-f_j)\|^2_\CH<\vare^2.
\end{equation*}
This yields that the vectors $(\Pi^\circ
\BG\phi_j)(E)\in\ell^2_\nu$ are linearly independent for $E$ lying
outside a subset $\CX_\vare$ of a small measure $\y(\vare)$, and
$\y(\vare)\to 0$ as $\vare\to 0$. For any
$E\in\D\setminus\CX_\vare$ we have $\rank \BV_{\BA,\BG}(E)=\nu$.
Letting $\vare\to 0$, we conclude from \eqref{Aa.bg} that
$\rank\BV_{\BA,\BG}=\nu$ a.e. on $\D$.

\vskip0.2cm

The equality \eqref{Aa.dim} is justified for any $\nu<\infty$.
Therefore, it remains valid also if $\nu=\infty$.

\bigskip{\bf Acknowledgments.}
The work on the paper started when S.N. visited the Weizmann
Institute in November -- December of 2003. The visit was supported
in part by the Department of Mathematics, and in part by the
Einstein center for theoretical physics. S.N. takes this
opportunity to express his gratitude to the Institute for its
hospitality and financial support. \vskip0.2cm M.S. acknowledges
partial financial support of the network SPECT of the ESF.

\bibliographystyle{amsplain}

\end{document}